\documentclass{amsart}
\usepackage{pinlabel}
\usepackage{amsmath, amsthm, amssymb, amscd, mathrsfs, eucal, epsfig,color}
\usepackage{hyperref}

\input xy
\xyoption{all}

\begin{document}

\newtheorem{theorem}{Theorem}[section]
\newtheorem{lemma}[theorem]{Lemma}
\newtheorem{corollary}[theorem]{Corollary}
\newtheorem{conjecture}[theorem]{Conjecture}
\newtheorem{proposition}[theorem]{Proposition}
\newtheorem{question}[theorem]{Question}
\newtheorem*{answer}{Answer}
\newtheorem{problem}[theorem]{Problem}
\newtheorem*{simplex_theorem}{Superideal Simplex Theorem}
\newtheorem*{cube_theorem}{Superideal Cube Theorem}
\newtheorem*{claim}{Claim}
\newtheorem*{criterion}{Criterion}
\theoremstyle{definition}
\newtheorem{definition}[theorem]{Definition}
\newtheorem{construction}[theorem]{Construction}
\newtheorem{notation}[theorem]{Notation}
\newtheorem{convention}[theorem]{Convention}
\newtheorem*{warning}{Warning}
\newtheorem*{assumption}{Simplifying Assumptions}

\theoremstyle{remark}
\newtheorem{remark}[theorem]{Remark}
\newtheorem{example}[theorem]{Example}
\newtheorem{scholium}[theorem]{Scholium}
\newtheorem*{case}{Case}


\def\ddal{\mathop{\vrule height 7pt depth0.2pt
\hbox{\vrule height 0.5pt depth0.2pt width 6.2pt}\vrule height 7pt depth0.2pt width0.8pt
\kern-7.4pt\hbox{\vrule height 7pt depth-6.7pt width 7.pt}}}
\def\sdal{\mathop{\kern0.1pt\vrule height 4.9pt depth0.15pt
\hbox{\vrule height 0.3pt depth0.15pt width 4.6pt}\vrule height 4.9pt depth0.15pt width0.7pt
\kern-5.7pt\hbox{\vrule height 4.9pt depth-4.7pt width 5.3pt}}}
\def\ssdal{\mathop{\kern0.1pt\vrule height 3.8pt depth0.1pt width0.2pt
\hbox{\vrule height 0.3pt depth0.1pt width 3.6pt}\vrule height 3.8pt depth0.1pt width0.5pt
\kern-4.4pt\hbox{\vrule height 4pt depth-3.9pt width 4.2pt}}}

\def\dal{\mathchoice{\ddal}{\ddal}{\sdal}{\ssdal}}

\def\id{\text{id}}
\def\Id{\text{Id}}
\def\1{{\bf{1}}}
\def\p{{\mathfrak{p}}}
\def\H{\mathbb H}
\def\Z{\mathbb Z}
\def\R{\mathbb R}
\def\C{\mathbb C}
\def\F{\mathbb F}
\def\P{\mathbb P}
\def\Q{\mathbb Q}
\def\E{{\mathcal E}}
\def\D{{\mathcal D}}
\def\d{{\mathfrak d}}
\def\A{{\mathcal A}}
\def\I{{\mathcal I}}
\def\density{\textnormal{density}}
\def\Square{\dal}

\def\tra{\textnormal{tr}}
\def\length{\textnormal{length}}

\def\cube{\textnormal{cube}}

\newcommand{\marginal}[1]{\marginpar{\tiny #1}}

\title{Coxeter groups and random groups}
\author{Danny Calegari}
\address{Department of Mathematics \\ University of Chicago \\
Chicago, Illinois, 60637}
\email{dannyc@math.uchicago.edu}
\date{\today}

\begin{abstract}
For every dimension $d$, there is an infinite family of convex cocompact reflection groups
of isometries of hyperbolic $d$-space --- the {\em superideal} (simplicial and cubical)
reflection groups --- with the property that
a random group at any density less than a half (or in the few relators model) contains
quasiconvex subgroups commensurable with some member of the family, with overwhelming probability.
\end{abstract}

\maketitle

\setcounter{tocdepth}{1}
\tableofcontents

\section{Introduction}

The recent combined work of Wise \cite{Wise}, Ollivier-Wise \cite{Ollivier_Wise} 
and Agol \cite{Agol} shows that random groups at density $< 1/6$ virtually embed in
right-angled Artin groups. Our main theorem is complementary to this (at least in flavor), 
and is concerned with virtually embedding certain classes of Coxeter groups in random groups (at any density).

\medskip

A {\em superideal} hyperbolic polyhedron is obtained by taking a polyhedron $P$ in real projective
space (of some dimension) and intersecting it with the region $H$ bounded by a quadric
(which is identified with hyperbolic space in the Klein model),
so that one obtains an infinite hyperbolic polyhedron $P\cap H$. If $P$ is a regular polyhedron
and $H$ is centered at the center of $P$, one obtains a {\em regular superideal polyhedron}
$P\cap H$, for which the symmetries of $P$ are realized as hyperbolic isometries of $P\cap H$.
For judiciously chosen $P$, the dihedral angles of $P\cap H$ might be of the form
$2\pi/m$ for some integer $m$, and hyperbolic space can then be tiled by copies of $P\cap H$.
The symmetry group of this tiling is a {\em superideal reflection group}.

In every dimension $d$, there are two interesting infinite families of convex cocompact
hyperbolic reflection groups --- the {\em superideal simplex} reflection group, and the
{\em superideal cube} reflection group. We denote the simplex groups by $\Delta(m,d)$ and
the cube groups by $\Square(m,d)$, where the dihedral angles of the superideal regular
polyhedra in each case are $2\pi/m$.

These reflection groups are Coxeter groups; for $\Delta(m,d)$ the Coxeter diagram is
\begin{center}
\begin{picture}(190,20)(0,0)
\put(5,10){\circle*{3}}
\put(35,10){\circle*{3}}
\put(65,10){\circle*{3}}
\put(95,10){\circle*{3}}
\put(155,10){\circle*{3}}
\put(185,10){\circle*{3}}
\put(5,10){\line(1,0){90}}
\multiput(95,10)(6,0){10}{\line(1,0){3}}
\put(155,10){\line(1,0){30}}
\put(16,15){$m$}
\end{picture}
\end{center}
whereas for $\Square(m,d)$ the Coxeter diagram is
\begin{center}
\begin{picture}(190,20)(0,0)
\put(5,10){\circle*{3}}
\put(35,10){\circle*{3}}
\put(65,10){\circle*{3}}
\put(95,10){\circle*{3}}
\put(155,10){\circle*{3}}
\put(185,10){\circle*{3}}
\put(5,10){\line(1,0){90}}
\multiput(95,10)(6,0){10}{\line(1,0){3}}
\put(155,10){\line(1,0){30}}
\put(16,15){$m$}
\put(166,15){$4$}
\end{picture}
\end{center}
where in each case there are $d+1$ nodes.

The main result we prove in this paper is that for any $d$, a random finitely presented
group at any density less than a half (or in the few relators model) contains quasiconvex
subgroups which are commensurable with some $\Delta(m,d)$ for some large $m$, and with
some $\Square(m,d)$ for some large (possibly different) $m$, with overwhelming probability.

The two cases are treated in a similar way, but the combinatorial details of the argument
are different in each case. The precise statements of our main theorems are the
Superideal Simplex Theorem and the Superideal Cube
Theorem, stated in \S~\ref{section:statement}. The proof of these theorems is carried out
in \S~\ref{section:spine} and \S~\ref{section:construction}. 

\medskip

These theorems generalize previous joint work of the author with Alden Walker \cite{Calegari_Walker}
and with Henry Wilton \cite{Calegari_Wilton}, and the architecture of the proof greatly resembles the
proofs of the main theorems in those papers. The new ideas in this paper are mainly combinatorial
in nature.

\section{Statements of the main theorems}\label{section:statement}

We first discuss the superideal simplex groups $\Delta(m,d)$.

\begin{example}
Some low-dimensional examples are familiar:
\begin{itemize}
\item{In any dimension $d$, taking $m=3$ gives the finite Coxeter group $A_{d+1}$, which is just
the symmetric group $S_{d+2}$. Similarly, in any dimension $d$, taking $m=4$ gives the
finite Coxeter group $BC_{d+1}$, the symmetry group of the $(d+1)$-dimensional cube
(or equivalently of the $(d+1)$-dimensional cross-polytope), known as the {\em hyperoctahedral group}.}
\item{In dimension $d=2$ there is an (ordinary) hyperbolic simplex with angles $2\pi/m$ whenever
$m\ge 7$. The groups $\Delta(m,2)$ are all commensurable, and are commensurable with the 
fundamental groups of all closed surfaces of genus at least 2.}
\item{In dimension $d=3$ the ideal regular simplex has $m=6$; the group $\Delta(6,3)$ is
commensurable with the fundamental group of the complement of the figure 8 knot. For
$m\ge 7$ the superideal ``simplices'' have infinite volume, but the groups they generate are
convex cocompact, with limit set a Sierpinski carpet. As $m \to \infty$ the limit sets converge
to the Apollonian gasket, with Hausdorff dimension approximately $1.3057$.}
\item{In dimension $d\ge 4$ the simplices are genuinely superideal whenever $m\ge 5$.}
\end{itemize}
\end{example}

Our first main theorem says that for every fixed dimension $d$, a ``generic'' finitely
presented group will contain some subgroup isomorphic to a finite index subgroup of some
element of the $\Delta(m,d)$ family. Explicitly, we show:

\begin{simplex_theorem}
For any fixed $d$, a random group at any density $<1/2$ or in the few relators model contains 
(with overwhelming probability) a subgroup commensurable with the Coxeter group $\Delta(m,d)$
for some $m\ge 7$, where $\Delta(m,d)$ is the superideal simplex group with Coxeter diagram
\begin{center}
\begin{picture}(190,20)(0,0)
\put(5,10){\circle*{3}}
\put(35,10){\circle*{3}}
\put(65,10){\circle*{3}}
\put(95,10){\circle*{3}}
\put(155,10){\circle*{3}}
\put(185,10){\circle*{3}}
\put(5,10){\line(1,0){90}}
\multiput(95,10)(6,0){10}{\line(1,0){3}}
\put(155,10){\line(1,0){30}}
\put(16,15){$m$}
\end{picture}
\end{center}
where there are $d+1$ nodes.
\end{simplex_theorem}

\begin{example}
Some special cases of this theorem were already known:
\begin{itemize}
\item{Taking $d=2$, this theorem implies that a random group contains surface subgroups (with overwhelming
probability). This is the main theorem of Calegari-Walker in the paper \cite{Calegari_Walker}.}
\item{Taking $d=3$, this theorem implies that a random group contains a subgroup isomorphic to the
fundamental group of a hyperbolic 3-manifold with totally geodesic boundary. This is the main
theorem of Calegari-Wilton in the paper \cite{Calegari_Wilton}. In fact, the Commensurability
Theorem proved in that paper, is {\em exactly} the statement of our main theorem in the
case $d=3$.}
\end{itemize}
\end{example}

We now discuss the superideal cube groups $\Square(m,d)$.

\begin{example}
The low-dimensonal examples are again familiar:
\begin{enumerate}
\item{In any dimension $d$, taking $m=3$ gives the finite coxeter group $BC_{d+1}$, the
hyperoctahedral group. Taking $m=4$ gives the {\em hypercubic honeycomb} $\tilde{C}_{d+1}$,
which is virtually abelian, and acts cocompactly on Euclidean space.}
\item{In dimension $d=2$ the group $\Square(m,2)$ is the symmetry group of the tiling by
regular right-angled $m$-gons, whenever $m\ge 5$. These groups are all commensurable with each
other, and with the fundamental groups of all closed surfaces of genus at least $2$ (and with
all $\Delta(m',2)$).}
\item{In dimension $d=3$ the group $\Square(5,3)$ is the symmetry group of the tiling by
regular right-angled dodecahedra; this is commensurable with the
(orbifold) fundamental group of the orbifold with underlying space the 3-sphere,
and cone angle $\pi$ singularities along the components of the Borromean rings, as made
famous by Thurston. For $m\ge 6$ the groups $\Square(m,3)$ are convex cocompact
but not cocompact. Similarly, in dimension $d=4$ the group $\Square(5,4)$ is the symmetry
group of the tiling by regular right-angled $120$-cells, whereas $\Square(m,4)$ is convex
cocompact but not cocompact for $m\ge 6$.}
\item{In dimension $d\ge 5$ the groups $\Square(m,d)$ are never cocompact when $m\ge 5$.}
\end{enumerate}
\end{example}

The analog of the Superideal Simplex Theorem for such groups is the following:

\begin{cube_theorem}
For any fixed $d$, a random group at any density $<1/2$ or in the few relators model contains 
(with overwhelming probability) a subgroup commensurable with the Coxeter group $\Square(m,d)$
for some $m\ge 5$, where $\Square(m,d)$ is the superideal cube group with Coxeter diagram
\begin{center}
\begin{picture}(190,20)(0,0)
\put(5,10){\circle*{3}}
\put(35,10){\circle*{3}}
\put(65,10){\circle*{3}}
\put(95,10){\circle*{3}}
\put(155,10){\circle*{3}}
\put(185,10){\circle*{3}}
\put(5,10){\line(1,0){90}}
\multiput(95,10)(6,0){10}{\line(1,0){3}}
\put(155,10){\line(1,0){30}}
\put(16,15){$m$}
\put(166,15){$4$}
\end{picture}
\end{center}
where there are $d+1$ nodes.
\end{cube_theorem}

The case $d=2$ reduces to the existence of closed surface subgroups, proved by
Calegari-Walker in \cite{Calegari_Walker}, but all other cases are new.

\medskip

We prove these theorems under the following two simplifying assumptions:
\begin{enumerate}
\item{that the length $n$ of the relations is divisible 
by a finite list of specific integers (implicitly depending on
the density $D$); and}
\item{that the number $k$ of free
generators is sufficiently large depending on $d$ in either case; for the 
Superideal Simplex Theorem we assume $2k-1\ge d+1$, whereas for the the Superideal
Cube Theorem we assume $2k-1\ge 2d+1$.}
\end{enumerate}
These assumptions simplify the combinatorics, but
they are superfluous and at the end of \S~\ref{subsection:simplifying} 
we indicate why they can be dispensed with,
and the theorems are true in full generality.

\section{Spines}\label{section:spine}

The proof of our main theorems is entirely direct and constructive; it depends on building
certain 2-dimensional complexes subject to local and global combinatorial constraints from
pieces which correspond to relators in the random group (in a precise way). In this
section we describe these combinatorial constraints.

\subsection{Sets as modules}

We frequently need to discuss finite sets together with an action of a finite group.
If a finite group $H$ acts on a set $A$ we say that $A$ is an {\em $H$-module}. The cases
of interest in this paper are:

\begin{itemize}
\item{the symmetric group $S_d$ acting on a $d$-element set by the standard permutation action; and}
\item{the hyperoctahedral group $BC_d$ acting on a $2d$-element set of $d$ pairs 
by the standard action, which permutes the pairs and acts on each pair as $\Z/2\Z$.}
\end{itemize}
By abuse of notation, we refer to a set together with such structure as an $S_d$-module 
or $BC_d$-module respectively. If $A$ is a $2d$ element set with a $BC_d$-module structure,
and $a$ is an element of $A$, the other element of the pair containing $a$ is said to be
{\em antipodal} to $a$.

\subsection{Graphs and 2-complexes}

Let $G$ be a random group with $k$ generators at density $D<1/2$ and length $n$. That means
$G$ is the group defined by a presentation
$$G:=\langle x_1,\cdots,x_k\; | \; r_1,r_2,\cdots,r_\ell\rangle$$
where $\ell = (2k-1)^{nD}$, and where each $r_i$ is chosen randomly and independently from the
set of all (cyclically) reduced words in the free group $F_k:=\langle x_1,\cdots,x_k\rangle$
of length $n$, with the uniform distribution. For an introduction to random groups see
\cite{Gromov} or \cite{Ollivier}.

We let $r=r_1$, and let $G_r$ denote the random 1-relator group with presentation
$$G_r:=\langle x_1,\cdots,x_k\; | \; r\rangle$$

The group $G$ is the fundamental group of a 2-complex $K$, with one vertex, with one edge for
each generator, and with one 2-cell for each relator. We denote the 1-skeleton of $K$ by $X$;
thus $X$ is a $k$-fold rose. The group $G_r$ is the fundamental group of a subcomplex
$K_r$, whose 1-skeleton is equal to $X$, and which is contains only one 2-cell, the cell
attached by the relation $r=r_1$.

\begin{definition}
A {\em graph over $X$} is a graph $\Sigma$ together with a simplicial immersion from
$\Sigma$ into $X$. 

A {\em morphism over $X$} between oriented graphs $\Sigma\to X$, $\Sigma'\to X$ over $X$ is a simplicial
immersion $\Sigma \to \Sigma'$ such that the composition $\Sigma \to \Sigma' \to X$ is
the given immersion from $\Sigma \to X$.
\end{definition}

If $\Sigma$ is a graph over $X$, the oriented edges of $\Sigma$ can be labeled by the $x_i$
and their inverses, by pulling back the labels from the edges of $X$. Conversely, a graph
with oriented edges labeled by the $x_i$ and their inverses, has the structure of a graph
over $X$ providing no adjacent oriented edges have inverse labels.

Let $L$ be a finite union of simplicial circles (i.e.\/ circles subdivided into edges) with 
oriented edges labeled by generators or their inverses, in such a way that the cyclic word
on each component of $L$ is $r$. Then $L$ is a graph over $X$. 

We will build a simplicial graph $\Sigma$ over $X$ together with a morphism 
$L \to \Sigma$ over $X$ satisfying certain conditions. These conditions are different for
the Superideal Simplex Theorem and for the Superideal Cube Theorem, but they have common
features, which we now describe.

\begin{figure}[htpb]
\labellist
\small\hair 2pt
\endlabellist
\centering
\includegraphics[scale=0.5]{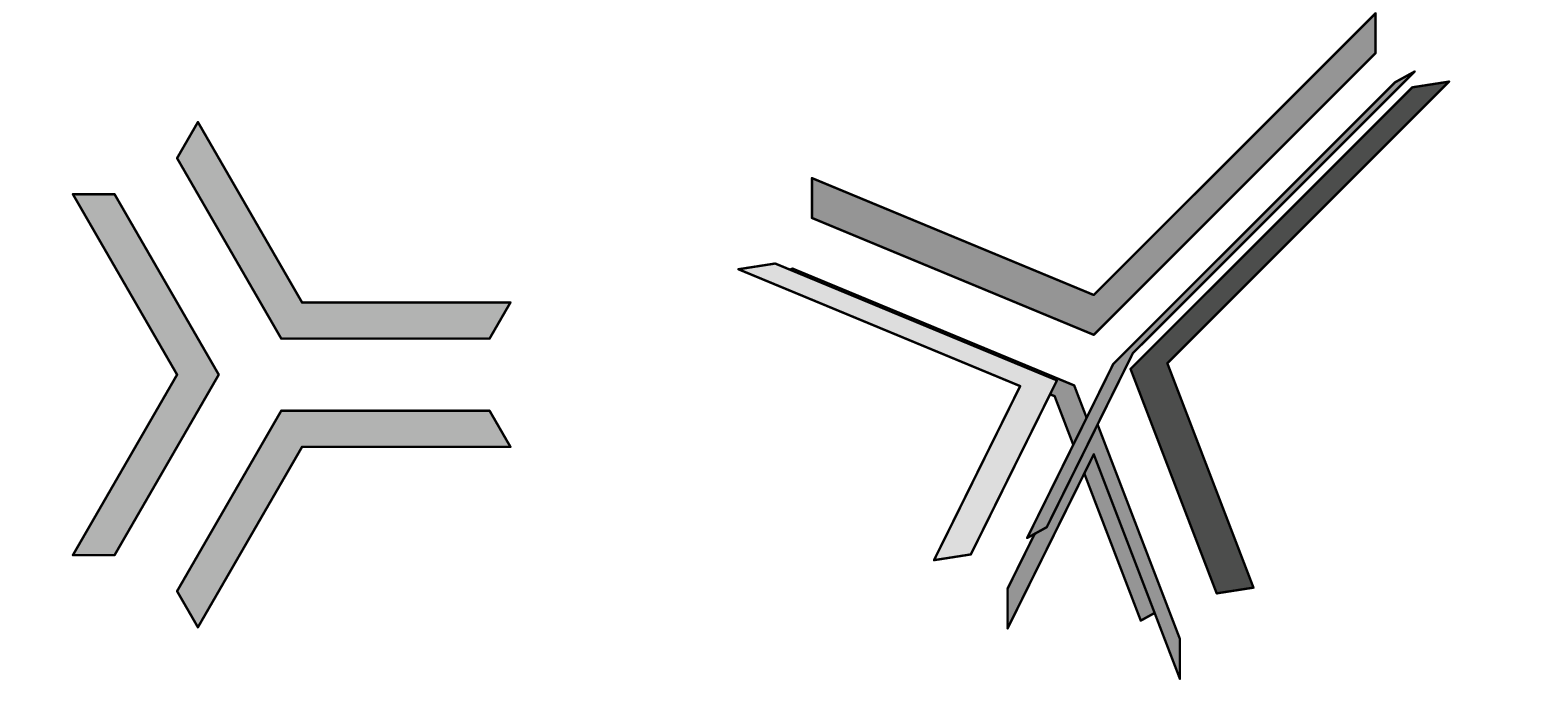}
\caption{local structure of simplicial spine near vertex in dimensions 2 and 3}\label{local_spine_simplex}
\end{figure}

\begin{figure}[htpb]
\labellist
\small\hair 2pt
\endlabellist
\centering
\includegraphics[scale=0.5]{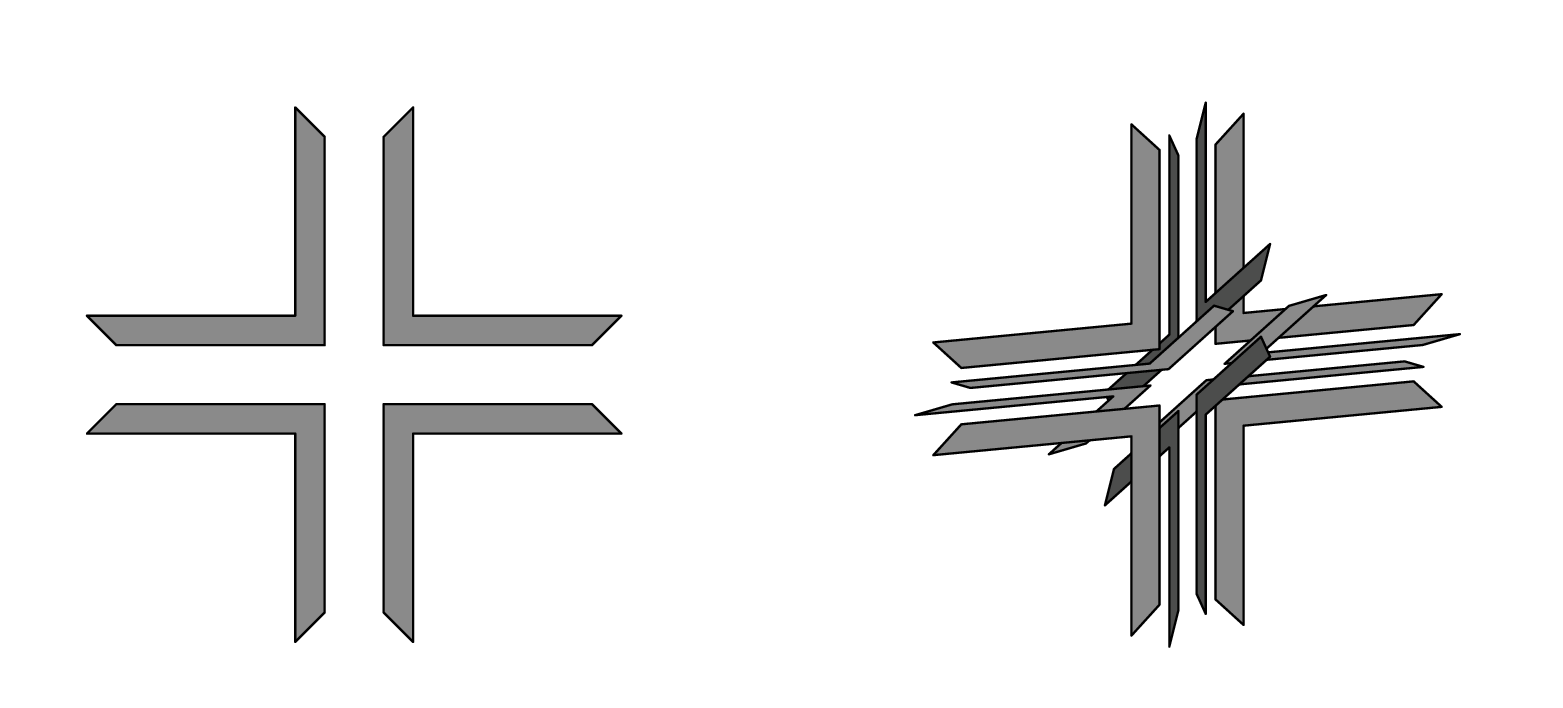}
\caption{local structure of cubical spine near vertex in dimensions 2 and 3}\label{local_spine_cube}
\end{figure}

\begin{definition}
An {\em $(m,d)$-regular simplicial spine} or {\em $(m,d)$-regular cubical spine} is a graph
$\Sigma$ over $X$, together with a morphism $L \to \Sigma$ over $X$ with the following
properties:
\begin{enumerate}
\item{The graph $\Sigma$ is obtained by subdividing a regular graph. That is, all the vertices
are either 2-valent (we call these {\em internal}) or of the same fixed valence $>2$
(we call these {\em genuine}). In the simplex case, each genuine 
vertex is $d+1$-valent; in the cube case, each genuine vertex is $2d$-valent 
and further admits the structure of a $BC_d$-module. By abuse of notation, we omit the
adjective ``genuine'' when discussing higher valence vertices unless the meaning would
be ambiguous.}
\item{The map $L \to \Sigma$ is $d$ to $1$ in the simplex case, and $2(d-1)$ to $1$ in
the cube case; i.e.\/ each edge of $\Sigma$ is the preimage of exactly $d$ (resp. $2(d-1)$)
edges of $L$. Furthermore, the set of preimages of each edge has the structure of
an $S_d$ (resp. $BC_{d-1}$)-module.}
\item{For each vertex $v$ of $\Sigma$ and each pair of distinct incident edges $e,e'$ 
which are not antipodal in the cube case, there is exactly one segment of $L$ of length $2$
whose midpoint maps to $v$ and whose adjacent edges map to $e$ and $e'$. 
In the cubical case, if $L(e)$ is the set of $2(d-1)$ edges of $L$ mapping to $e$ with its
$BC_{d-1}$-module structure, then the $2(d-1)$ edges in $L$ adjacent to $L(e)$ map 
bijectively to the $2(d-1)$ edges of $\Sigma$ not equal to $e$ or its antipode; this latter set also has a
natural $BC_{d-1}$-module structure by restriction, and we require that this map of
$BC_{d-1}$-modules respect the module structure; see Figure~\ref{local_spine_simplex} for the
local picture in the simplicial case and Figure~\ref{local_spine_cube} for the cubical case in
dimensions 2 and 3.}
\item{Each component of $L$ has exactly $m$ vertices which map to genuine vertices of 
$\Sigma$.}
\item{$\Sigma$ satisfies the {\em cocycle condition} (to be defined shortly).}
\end{enumerate}
\end{definition}

Let $L \to \Sigma$ be a spine satisfying the first three conditions to be $(m,d)$-regular. 
The mapping cylinder $M$ is a 2-complex which admits a canonical
local embedding into $\R^d$ in such a way that the (combinatorial) symmetries of the link of
each vertex or edge are realized by isometries of the ambient embedding (taking
module structure into account). Taking a tubular
neighborhood of this canonical thickening gives rise to a $D^{d-2}$-bundle over $M$ with
a flat connection with holonomy in the symmetric group $S_{d-1}$ in the simplicial case,
or $BC_{d-2}$ in the cubical case, acting by the standard
representation. The restriction of this bundle to each component of $L$ therefore determines
a conjugacy class in $S_{d-1}$ or $BC_{d-2}$ respectively. 
The {\em cocycle condition} is the condition that this
conjugacy class is trivial, for each component of $L$.

Let $L\to \Sigma$ be an $(m,d)$-regular spine. Form the mapping cylinder $M$ and then let $\overline{M}$ be
the 2-complex obtained from $M$ by attaching a disk to each component of $L$. Since each
component of $L$ has the label $r$, the
tautological immersion $L \to X$ extends to a cellular immersion $\overline{M} \to K_r$.

The defining properties of an $(m,d)$-regular simplicial (resp. cubical) spine ensure that
the 2-complex $\overline{M}$ is locally isomorphic to the 2-dimensional spine $\mathcal{S}$
which is dual to the tiling of $d$-dimensional hyperbolic space by regular 
superideal simplices (resp. cubes) with dihedral angles $\pi/m$. 
To see this, observe that the vertex stabilizers in
the symmetry group of $\mathcal{S}$ are isomorphic to $S_d$ (resp. $BC_{d-1}$) and their
action on vertex links of $\mathcal{S}$ agrees with the module structure on the 
vertex links of $\overline{M}$. There is therefore a developing map from the universal
cover of $\overline{M}$ to $\mathcal{S}$ which is an immersion, 
and therefore an isomorphism. This realizes the fundamental group of $\overline{M}$ as a
subgroup of $\Delta(m,d)$ (resp. $\Square(m,d)$) which is of finite index, since $\overline{M}$
is compact.

We will show in the next section that for $r$ a sufficiently
long random word, we can construct an $(m,d)$-regular spine (of either kind) in which all the 
maximal paths in $\Sigma$ containing only internal vertices in their interior
(informally, the ``topological edges'' of $\Sigma$) 
are bigger than any specified constant $\lambda$. This is the main
ingredient necessary to show that the immersion $\overline{M} \to K_r \to K$ is $\pi_1$-injective,
and stays injective when a further $(2k-1)^{nD}-1$ random relators are added.

\section{Constructing spines}\label{section:construction}

We start with $L$, a finite union of simplicial circles each labeled with the relator $r$.
The spine $\Sigma$ is obtained from $L$ by identifying oriented edges of $L$ in groups of $d$
or $2(d-1)$, respecting labels, in such a way as to satisfy the conditions of a 
regular spine. We
thus obtain $\Sigma$ as the limit of a (finite) sequence of quotients
$L \to \Sigma_i$ limiting to $\Sigma$, where each $L \to \Sigma_i$
is a morphism over $X$. Thus, each intermediate $\Sigma_i$ should immerse in $X$. We call
an identification of subgraphs of $\Sigma_{i-1}$ with the same labels, giving rise to a quotient
$\Sigma_{i-1} \to \Sigma_i$ over $X$ a {\em legal quotient} if $\Sigma_i$ is immersed in $X$.
We assume in the sequel that our quotients are always legal.

In our intermediate quotients $\Sigma_i$, some edges are in the image of $d$ (resp. $2(d-1)$) 
distinct edges of $L$, and some are in the image of a unique edge of $L$. We call edges of 
the first kind {\em glued}, and edges of the second kind {\em free}.

In order to simplify notation, in the sequel we let $\d$ denote either $d$ in the simplicial
case, or $2(d-1)$ in the cubical case.

We choose a big number $\lambda$ divisible by $\d$, and an even number $N$ with $N\gg \lambda$, 
and assume for the sake of convenience that $\lambda\cdot N$ divides
$n$, the length of the relator $r$ (the size of $\lambda$ that we need will ultimately depend only
on the density $D$, but not on the length $n$ of the relators).
Then we subdivide each component of $L$ into consecutive
{\em segments} of length $\lambda$.

\subsection{Creating beachballs}

A {\em block} in $L$ is a sequence of $N$ consecutive segments. Since we assume that $\lambda\cdot N$
divides $n$, each component of $L$ can be subdivided into $n/(\lambda\cdot N)$ consecutive blocks.
Each block is made up of alternating {\em odd segments} and {\em even segments}; by abuse
of notation we refer to the segment immediately before
the first segment in a block (which is contained in the previous block)
as the ``first even segment''.
We say that a $\d$-tuple of blocks is {\em compatible} if their odd segments can all be identified
in groups of $\d$ (in the order in which they appear in the blocks)
in such a way that the resulting quotient is legal. This means that each $\d$-tuple of
even segments at the same location in the $\d$ blocks
must start and end with {\em different} letters, and the same must be true for the last
$\d$-tuple of even segments in the blocks immediately preceding the given collection of blocks.
The existence of a compatible collection of blocks requires $2k-1\ge \d$. 

\begin{figure}[htpb]
\labellist
\small\hair 2pt
\endlabellist
\centering
\includegraphics[scale=0.75]{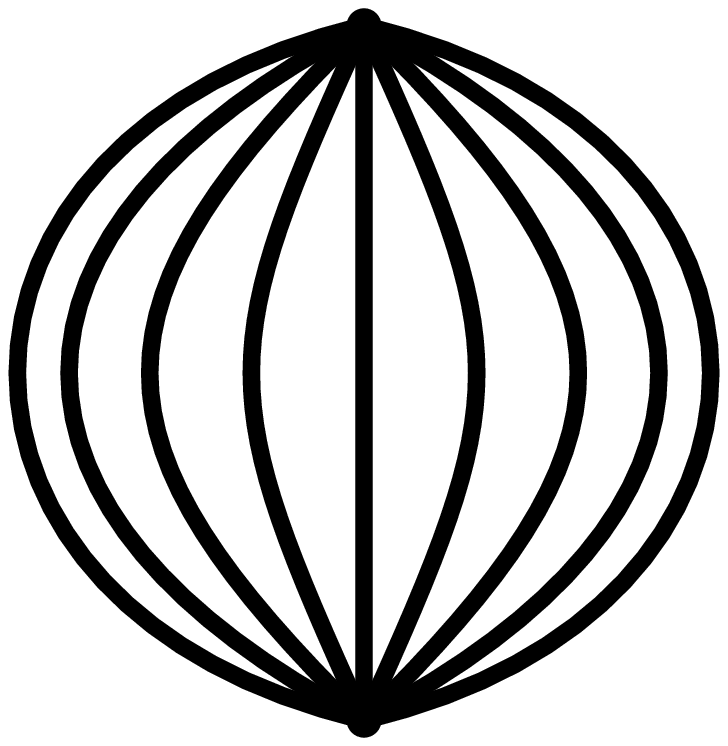}
\caption{A beachball for $\d=9$}\label{beachball}
\end{figure}

When $\d$ blocks are glued together along their odd segments, each collection of even segments
{\em except} for the first and last one is identified at their vertices, producing a {\em beachball} ---
a graph with two vertices and $\d$ edges, each joining one vertex to the other. See 
Figure~\ref{beachball}.

We say that a $(\d+1)$-tuple of blocks is {\em supercompatible} if each sub $\d$-tuple is 
compatible. Evidently, the existence of supercompatible tuples requires $2k-1\ge \d+1$, which
is our second simplifying assumption.

Recall that $\d=d$ in the simplicial case, and $\d=2(d-1)$ in the cubical case. 
Suppose $r$ is a random reduced word of length $n$. Since $r$ is random, if $n$ is sufficiently long
it is possible to partition a proportion $(1-\epsilon)$ 
of the blocks into compatible $\d$-tuples, and glue their odd segments, 
for $\epsilon$ depending only on $n$, and therefore
$1/\epsilon \gg N,\lambda$, with probability $1-O(e^{-n^c})$.
This leaves a proportion $\epsilon$ of the blocks left unglued. 
For each of these unglued blocks, find a
collection of $\d$ glued blocks so that the $\d+1$ tuple is supercompatible. If we unglue the
collection of $\d$ glued blocks, and then take $\d+1$ copies of our graph, we can glue the
resulting $\d(\d+1)$ blocks in compatible groups of $\d$.

In the cubical case, in order to give each collection of edges the natural structure
of a $BC_{d-1}$-module, it is convenient to first take some finite number of disjoint copies of
$L$ and partition them into $2(d-1)$ subsets of equal size. The group $BC_{d-1}$ acts on
each set of $2(d-1)$ copies of $L$ as a standard representation, giving these components
the structure of a $BC_{d-1}$-module. When we partition blocks into $\d$-tuples (where 
$\d=2(d-1)$ in the cubical case), each tuple should contain exactly
one block from each subset, so that the elements of the block inherit a natural $BC_{d-1}$
action.  

Thus at the end of this step, we have constructed $\Sigma_1$ which consists of collections of
glued edges of length $\lambda$ (each the image of $\d$ edges of $L$), 
a {\em reservoir} of beachballs with every edge of length
$\lambda$ and all possible edge labels in almost uniform distribution of possible edge labels, 
and a {\em remainder} consisting of a $\d$-valent graph in which every edge has length $\lambda$.
Note that the remainder is obtained as the union of the first and last (even) segments from each block,
so the mass of the remainder is of order $O(1/N)$, where by assumption we have
chosen $N\gg \lambda$. Furthermore in the cubical case, for each component $X$ 
of the reservoir or remainder, the edges of $X$ incident to each vertex inherit a natural
$BC_{d-1}$-module structure.

\subsection{Covering move}\label{subsection:covering}

The next move is called the {\em covering move}, and its effect is to {\em undo} some of
the gluing of some previous step, and then to {\em reglue} in such a way that the net effect is
to transform some collection of beachballs into more complicated graphs, each of which is
a finite cover of a beachball, with prescribed topology and edge labels. This move does
{\em not} necessarily preserve the total mass of beachballs, and can be used to adjust the
total mass of beachballs of any specific kind, at the cost of transforming some other
(prescribed) set of beachballs into covers.

\begin{figure}[htpb]
\labellist
\small\hair 2pt
\endlabellist
\centering
\includegraphics[scale=0.4]{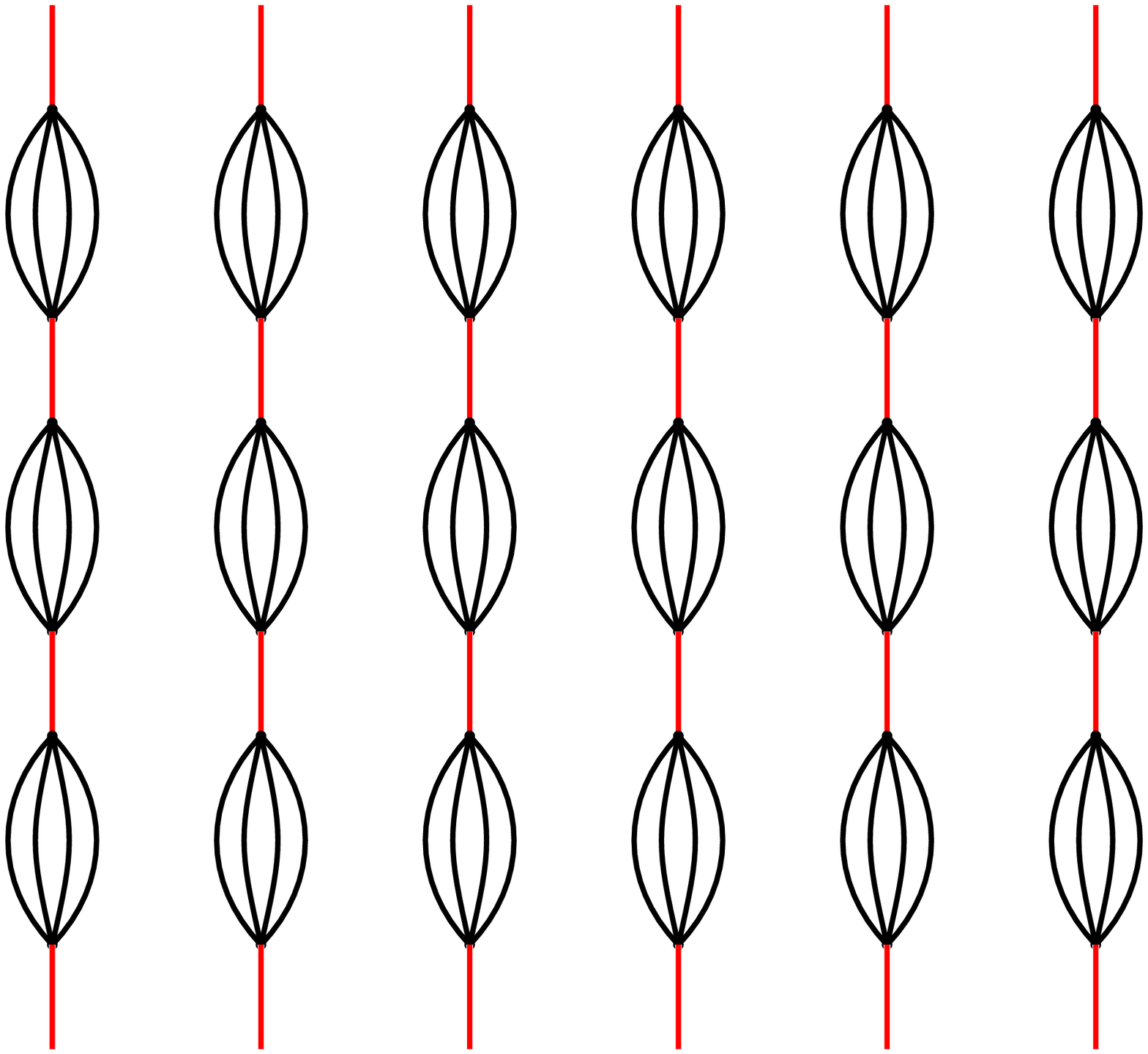}
\caption{A matrix of 7 columns of 3 consecutive beachballs. 
In this and subsequent figures, glued edges are in red and free edges are in black.}\label{matrix}
\end{figure}

To describe this move, we first take some consecutive strings of beachballs and arrange them
in a ``matrix'' form, where each string of beachballs is a column of the matrix, and the
different strings make up the rows of the matrix. See Figure~\ref{matrix} for a matrix
consisting of 7 columns of 3 consecutive beachballs.

The covering move takes as input a matrix of
$s$ columns of $r$ consecutive beachballs with compatible
labels on corresponding glued segments. This means the following. In each column there
are $r+1$ consecutive glued segments interspersed with the $r$ beachballs. Thus
between each pair of consecutive rows of beachballs there is a row of glued segments, each
consisting of $\d$ glued edges with some label. We require that all the labels on glued edges
in the same row are equal.

In the cubical case, the edges in each column have a natural structure of a
$BC_{d-1}$-module. Thus it is possible to ``trivialize'' the module structure on the
union of edges in each row, giving them the structure of a product of a standard
$BC_{d-1}$-module with an $s$-element set. 

The covering move pulls apart each interior row of glued segments (i.e.\/ the rows between
a pair of consecutive row of beachballs) into $\d s$ edges, all with the same label,
then permutes them into $s$ new sets of $\d$ elements and reglues them. More specifically,
if we identify the set of $\d s$ edges with the product $\lbrace 1,\cdots,s\rbrace \times 
\lbrace 1,\cdots,\d\rbrace$ where $BC_{d-1}$ acts on the second factor, we permute each
$\lbrace 1,\cdots,s\rbrace \times i$ factor by some element of $S_s$. Thus the move
is described on each of the rows by an element $p_j \in S_s^\d$, the product of $\d$
copies of the symmetric group $S_s$, for $0\le j \le r$ where
$p_0$ and $p_r$ are the identity permutation.

For each $1\le j \le r$ the $s$ beachballs in the $j$th row are rearranged by this move
into a graph which is a (possibly disconnected) covering space of a beachball of degree $s$;
which covering space is determined by the permutations $p_{j-1}$ and $p_j$. If
$p_{j-1} = p_j$ it is a trivial covering space, in the sense that it consists just of
$s$ disjoint beachballs. However, even if it is trivial as a covering space,
if $p_{j-1}$ and $p_j$ are not the identity, the
edges making up each beachball might be different from the edges making up the beachballs
before the move. We allow two possibilities for the labels on the new beachball covers
in each row:
\begin{enumerate}
\item{the beachball cover in the given row can have arbitrary topology, 
and the labels on the edges are legal; or}
\item{the beachball cover in the given row is trivial (i.e.\/ $p_{j-1}=p_j$) and on each
beachball the edges all have the {\em same} label.}
\end{enumerate}

In the former case we do nothing else to the row. In the latter case, since all the edge
labels are equal (and the gluing is not legal as it stands), 
we collapse the beachball to the interior of a glued segment of length $3\lambda$. So
the net effect of the covering move is to transform some set of beachballs into covering
spaces of the same mass, while possibly eliminating some set of beachballs of comparable
mass.

Two kinds of beachball covering spaces are especially important in what follows; these
are 

\begin{figure}[htpb]
\labellist
\small\hair 2pt
\endlabellist
\centering
\includegraphics[scale=0.4]{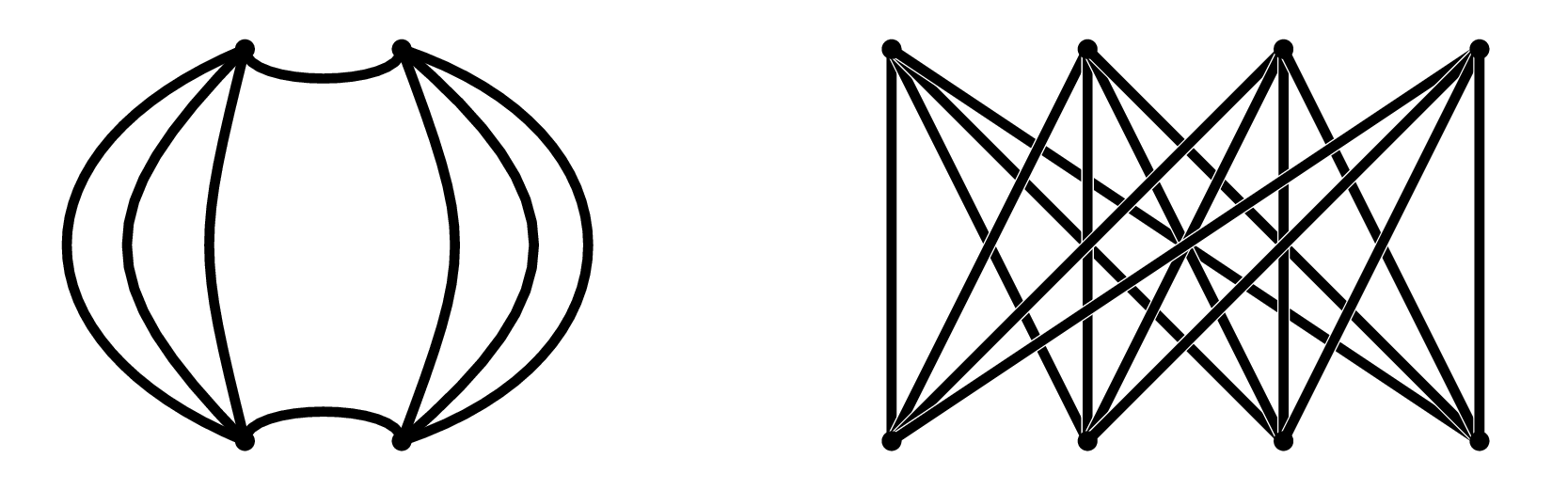}
\caption{A barrel and a bipart for $\d=4$}\label{barrel_bipart}
\end{figure}

\begin{itemize}
\item{a 2-fold cover of a beachball called a {\em barrel}; and}
\item{a $\d$-fold cover of a beachball called a {\em bipart}, whose underlying graph is
the complete bipartite graph $K_{\d,\d}$.}
\end{itemize}
See Figure~\ref{barrel_bipart}.

We now describe two especially useful examples of covering moves which will be used in the
sequel.

\begin{example}[Elimination]\label{example:elimination}
The elimination move trades $2\d$ beachballs for 2 biparts, and eliminates a collection of $\d$
beachballs with the same labels. It is convenient to perform this move in such a way that the
biparts created have labels of {\em covering type}. This means that the edge labels
on the bipart are pulled back from some legal edge labeling on a beachball under the
covering projection. See Figure~\ref{adjust}.

\begin{figure}[htpb]
\labellist
\small\hair 2pt
\endlabellist
\centering
\includegraphics[scale=0.29]{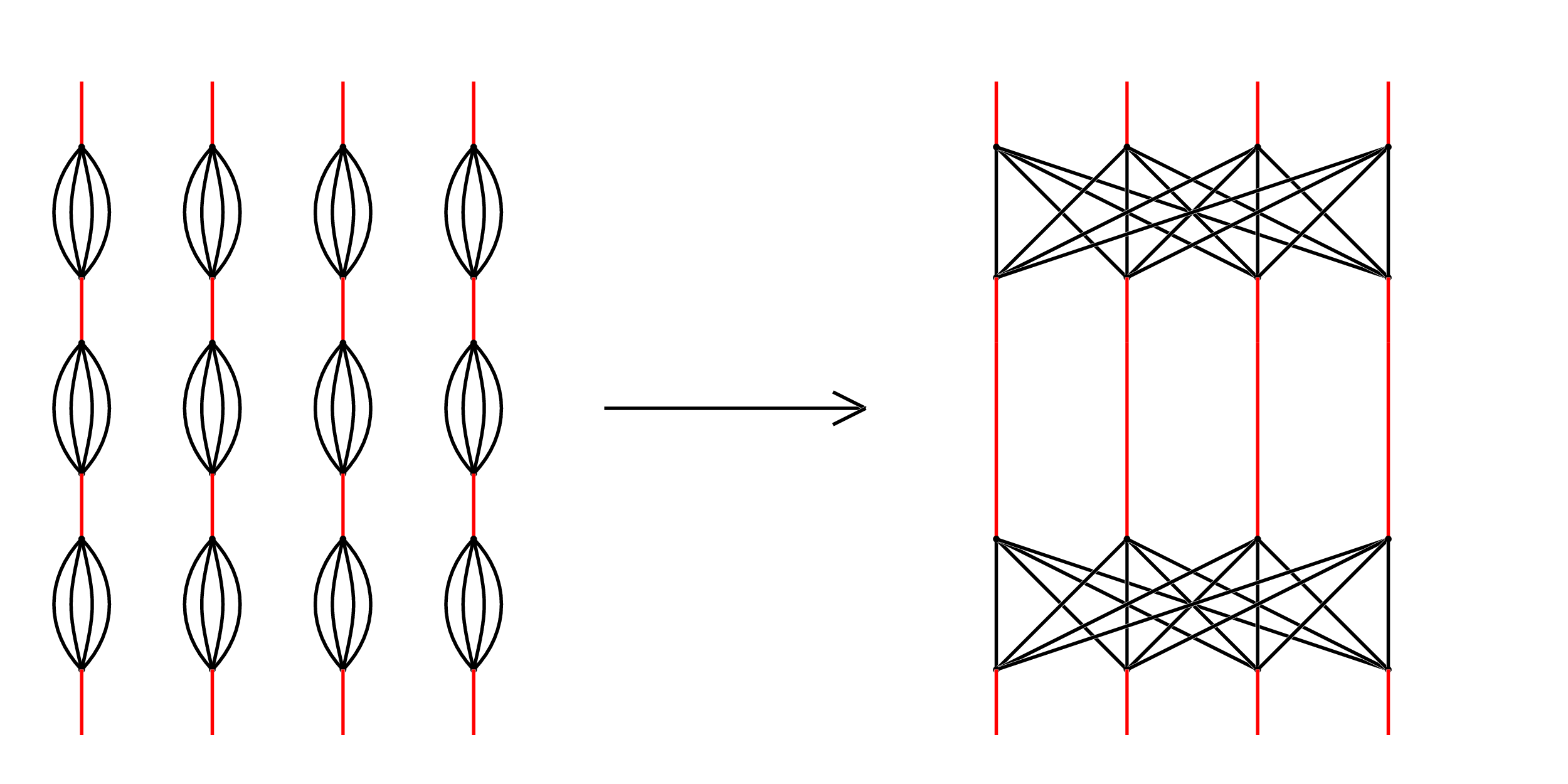}
\caption{The elimination move for $\d=4$}\label{adjust}
\end{figure}
\end{example}

\begin{example}[Rolling barrels]\label{example:rolling}
The rolling move trades 4 beachballs for 2 barrels. This move is self-explanatory.
See Figure~\ref{rolling}.

\begin{figure}[htpb]
\labellist
\small\hair 2pt
\endlabellist
\centering
\includegraphics[scale=0.4]{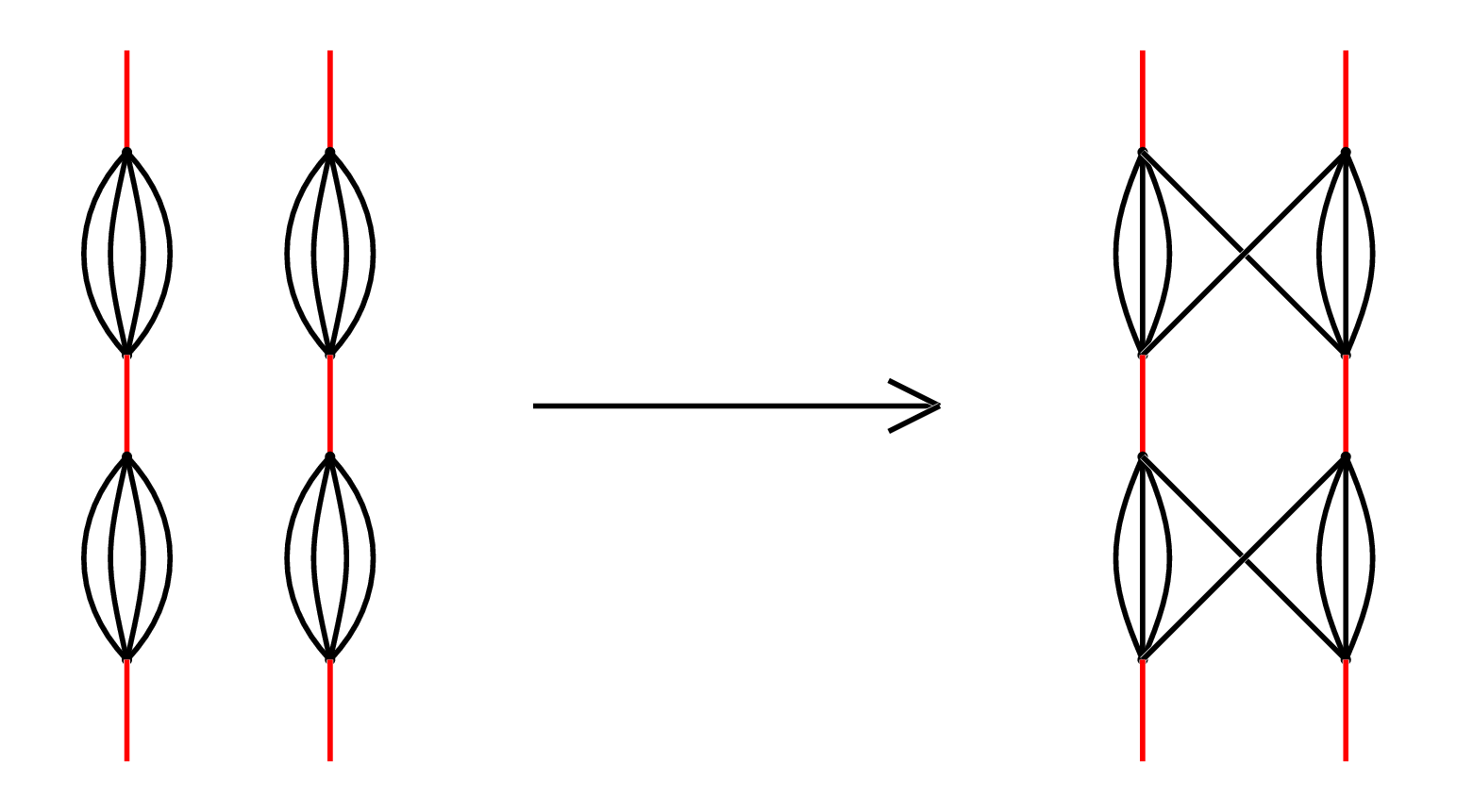}
\caption{The rolling move for $\d=4$}\label{rolling}
\end{figure}
\end{example}

\subsection{Tearing the remainder}\label{subsection:tear}

We now describe a modification of the gluing which replaces the remainder with a new remainder
of similar mass, but composed entirely of barrels. In the process,
a set of beachballs of similar mass are transformed into biparts.

\begin{figure}[htpb]
\labellist
\small\hair 2pt
\endlabellist
\centering
\includegraphics[scale=0.29]{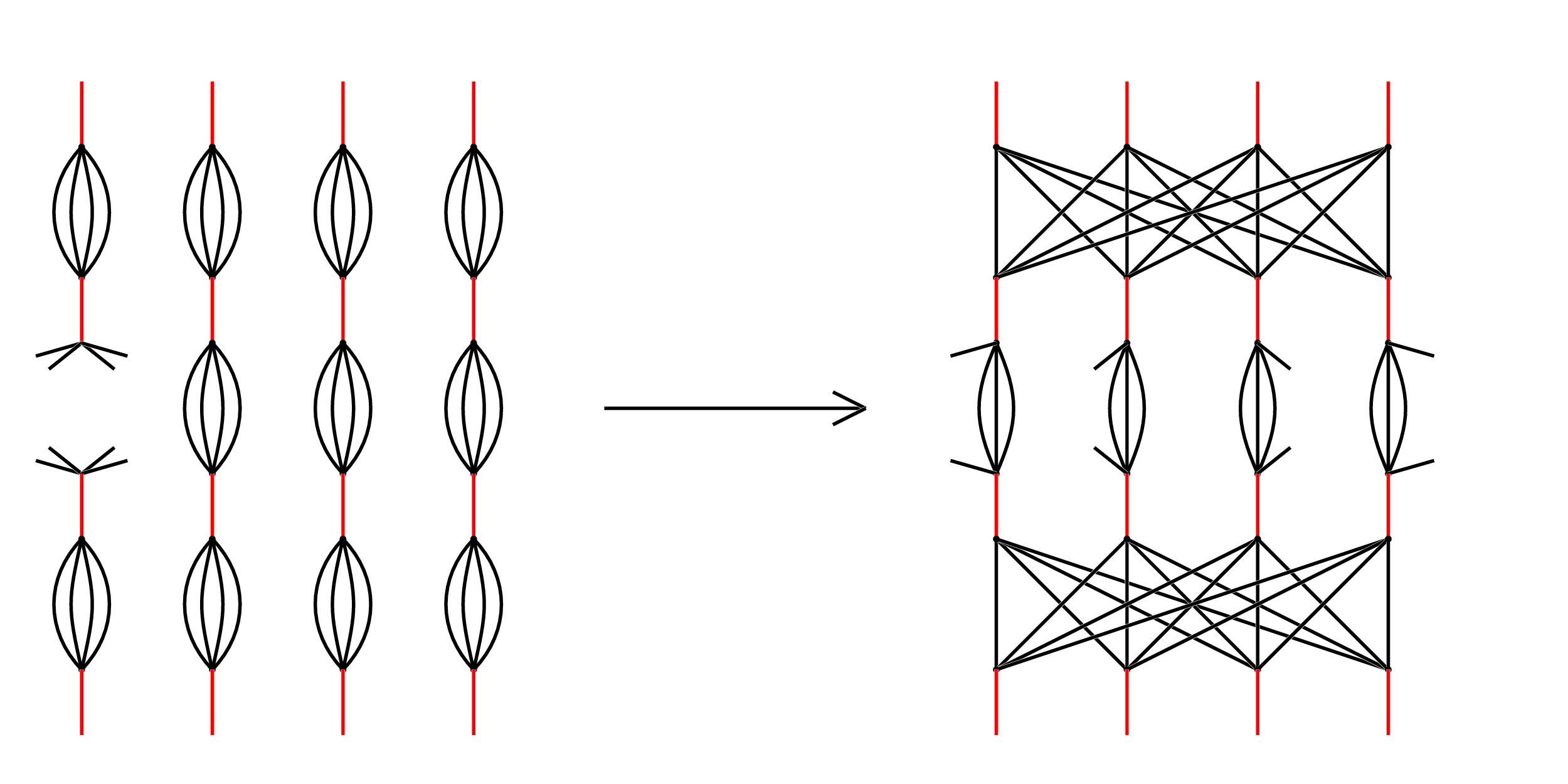}
\caption{The tear move for $\d=4$. }\label{tear}
\end{figure}

The tear move is applied to two vertices $v,v'$ in a component $Y$ of the remainder. It uses up
the two beachballs in the reservoir which are adjacent to $v$ and $v'$, and a further
$3(\d-1)$ beachballs drawn from the middle of the reservoir in consecutive
groups of three. The result of the move ``tears'' $Y$ apart at
$v$ and $v'$, producing $\d$ copies of each vertex $v_i,v_i'$ for $1\le i\le \d$, each
$v_i$ joined to $v_i'$ by ${\d-1}$ new edges (these new $\d(\d-1)$ edges coming from $(\d-1)$ of the
beachballs), and transforms the remaining $2\d$ beachballs into two biparts. In the process,
$2\d$ glued edges are unglued, and then reglued in a different configuration; thus it is necessary
for the labels on each set of $\d$ edges to agree.
The tear move is illustrated in Figure~\ref{tear}; as in \S~\ref{subsection:covering},
the segments before and after the move have the same vertical coordinate, and the gluing
respects $BC_{d-1}$-module structure on edges in the cubical case.

\begin{figure}[htpb]
\labellist
\small\hair 2pt
\endlabellist
\centering
\includegraphics[scale=0.4]{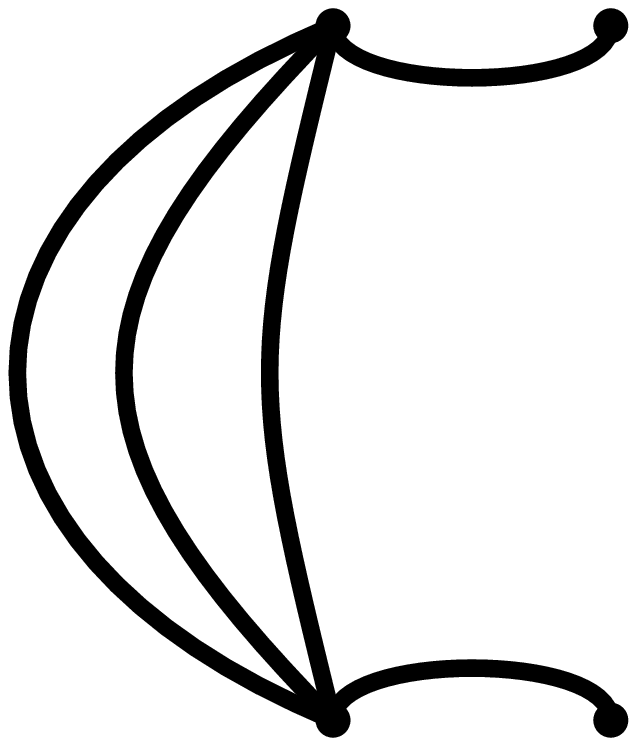}
\caption{A half-barrel for $\d=4$. }\label{half_barrel}
\end{figure}

Applying a tear move at a pair of vertices $v,v'$ in $Y$ creates a collection of {\em half-barrels}
i.e.\/ a beachball with one edge replaced by a pair of {\em free edges} joined to nearby
{\em free vertices}; see Figure~\ref{half_barrel}. 

We refer to the non-free edges in a half-barrel as {\em barrel edges}, and
observe that the number of non-barrel edges in $Y$ is unchanged by a tear move. Thus after
applying the tear move once to create several half-barrels, we can apply the tear move again
to the pair of free vertices at the end of one of the half-barrels, in the process creating an
honest barrel which may be transferred to the reservoir, reducing the total number of
non-barrel edges in $Y$. We can use total number of non-barrel edges of $Y$ as
a measure of complexity, and observe that after applying $O(1/N)$ tear moves, 
we can tranform the entire
remainder to a mass of $O(1/N)$ barrels and biparts, while the reservoir still contains
almost exactly the same number of beachballs in almost exactly the uniform distribution.

\subsection{Hypercube gluing}\label{subsection:hypercube_gluing}

At this stage of the construction we have a reservoir consisting of an almost equidistributed
collection of beachballs, and a relatively small mass of pieces consisting of biparts
and barrels. In this section and the next we explain how to glue up beachballs. There are
two different (but related) constructions depending on whether we are in the simplicial 
or cubical case. In the simplicial case we use the {\em hypercube gluing}.

A finite cover of the hypercube gluing can be used to glue up all the barrels and biparts 
produced by the tear moves, together with some other complementary collection of pieces 
produced by the covering move. The existence of such a combinatorial cover in which
the pieces from the remainder can be included follows by LERF for free groups, and 
it is straightforward to use the covering move to produce complementary pieces of the
correct topology and labels. An explicit description of the combinatorics of the
necessary covers in dimension 3 is given in \cite{Calegari_Wilton}, \S~3.10, and the
general case is very similar.

After dealing with the remainder, almost all the mass of the reservoir can be completely 
glued up using hypercube gluing, and the error can be transformed away using
the elimination move, while simultaneously transforming some prescribed collection
of beachballs into biparts of covering type, which can themselves be glued up
by a covering of the hypercube gluing.

\begin{figure}[htpb]
\labellist
\small\hair 2pt
\endlabellist
\centering
\includegraphics[scale=0.7]{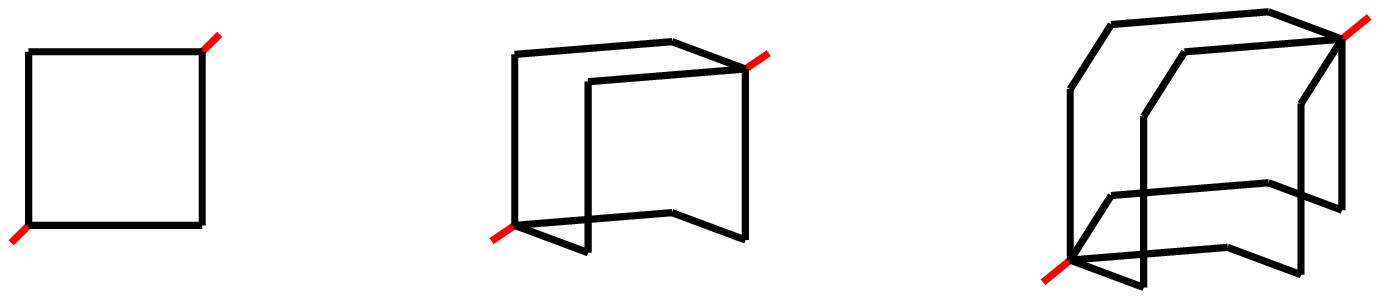}
\caption{Beachball immersed in the 1-skeleton of a cube in dimensions 2, 3 and 4.}\label{hypercube}
\end{figure}

It remains to describe the hypercube gluing (recall that we are in the simplicial case).
Let $C$ denote the $1$-skeleton of the $d$-dimensional cube, and with each edge
of length $\lambda/d$. 
Let $e_1,\cdots,e_d$ denote vectors of length $\lambda/d$ aligned positively along the $d$
coordinate axes in order. There is a path of length $\lambda$ from the vertex 
$(0,0,\cdots,0)$ to the vertex $(\lambda/d,\lambda/d,\cdots,\lambda/d)$ of $C$, obtained by
concatenating straight segments of length $\lambda/d$ in the order $e_1,e_2,\cdots,e_d$.
Taking $d$ cyclic permutations of this sequence of paths gives $d$ paths of length $\lambda$
between extremal vertices of $C$, whose union is a beachball. This is illustrated in
Figure~\ref{hypercube} in dimensions 2, 3 and 4.

There are $2^{d-1}$ pairs of
extremal vertices of $C$, and the union of $2^{d-1}$ suitably labeled beachballs can be arranged
along the $1$-skeleton of $C$ as above. Gluing the beachballs in this pattern is legal, and the
resulting graph locally satisfies the first three conditions to be a 
$(m,d)$-regular simplicial spine.

\subsection{Lens gluing}

The analog of the hypercube gluing for cubical spines is the {\em lens gluing}. This
is the most complicated move in our construction, and we first describe this move 
in low dimensions before giving the definition in generality.

\begin{figure}[htpb]
\labellist
\small\hair 2pt
\endlabellist
\centering
\includegraphics[scale=1]{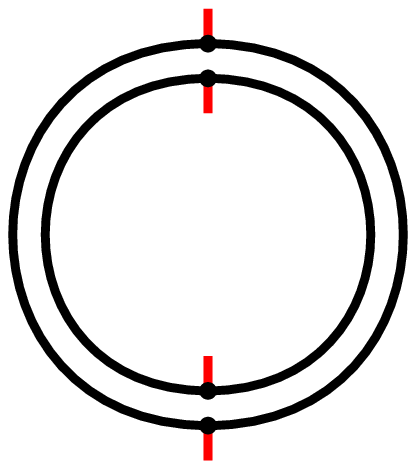}
\caption{Lens gluing $d=2$}\label{2lens}
\end{figure}

In dimension $2$ a beachball has degree 2 (i.e.\/ it has 2 edges); $BC_1=\Z/2\Z$ acts as the full
permutation group. Two beachballs with suitable labels can be identified along their
boundaries, laid out along a circle; see Figure~\ref{2lens}.

In dimension $3$ a beachball has degree 4; $BC_2=D_4$ acts by the dihedral group,
preserving or reversing a circular order on the 4 edges.
Eight beachballs with suitable labels can be glued up along a 
$K_{4,4}$ graph in $S^3$, where each edge of each beachball runs along a segment of length
$2$ in the $K_{4,4}$ graph, as in Figure~\ref{3lens} (where one of the vertices is at ``infinity''
in the figure):

\begin{figure}[htpb]
\labellist
\small\hair 2pt
\endlabellist
\centering
\includegraphics[scale=0.8]{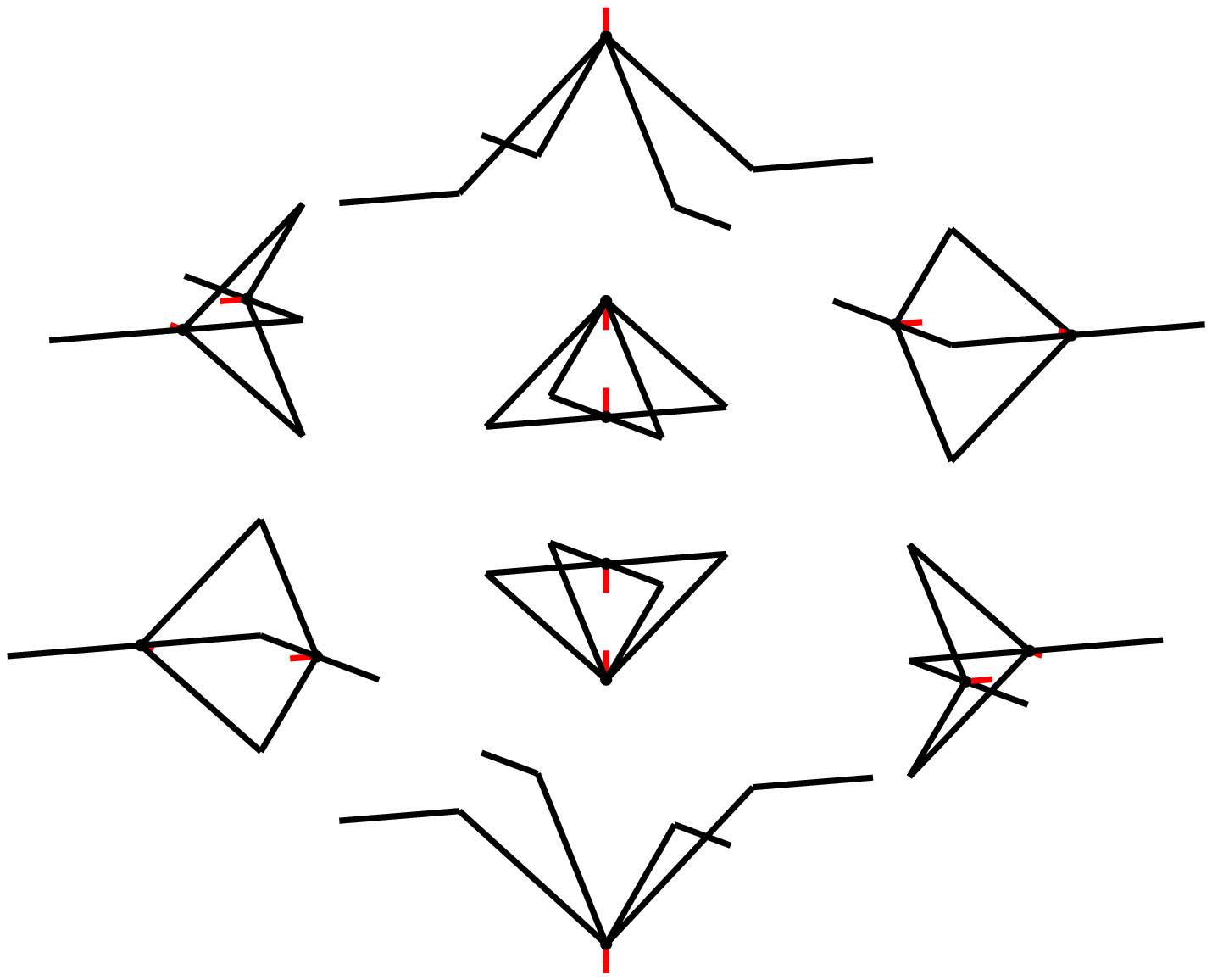}
\caption{Lens gluing $d=3$}\label{3lens}
\end{figure}

The lens gluing in dimension $d$ is closely related to the hypercube gluing in dimension $d-1$.
Recall the definition of the hypercube gluing move from \S~\ref{subsection:hypercube_gluing}.
In this move, $2^{d-2}$ beachballs of degree $d-1$ are draped along the 1-skeleton of a $d-1$-dimensional
cube; each edge of each beachball is subdivided into $d-1$ segments which run between antipodal
vertices of the cube. We let $C$ denote the 1-skeleton of the cube, and denote the immersion
of a beachball $b$ into $C$ by $b \to C$.

For any graph $L$, the {\em spherical graph} associated to $L$ is the graph with one $S^0$
(i.e.\/ two vertices) for each vertex in $L$, and one $S^0 * S^0$ for each edge in $L$. If
$v$ is a vertex of $L$, we denote the corresponding pair of vertices in $S(L)$ by $v0$ and $v1$.
Similarly, if $e$ is an edge of $L$, we denote the corresponding edges in $S(L)$ by
$e00$, $e01$, $e10$ and $e11$ (in fact this construction generalizes in an obvious way to
simplicial complexes of arbitrary dimension). There is a canonical simplicial map
from $S(L) \to L$ which forgets the $01$ labels. 

Now we apply this construction to $C$; the lens gluing will map $4\times 2^{d-2}$ beachballs
of degree $2(d-1)$ into $S(C)$ in such a way that the projection to $C$ will map each 4
beachballs onto each degree $(d-1)$ beachball mapping into $C$ in the hypercube gluing, with
each edge of the beachball downstairs in the image of two edges upstairs.

\begin{figure}[htpb]
\labellist
\small\hair 2pt
\endlabellist
\centering
\includegraphics[scale=0.7]{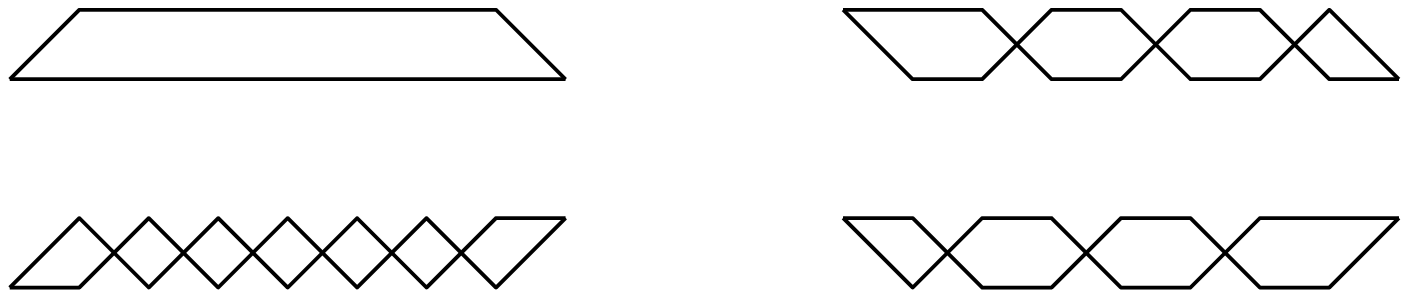}
\caption{Height graphs for $d=9$}\label{height_graph}
\end{figure}

If $b \to C$ is a beachball, and $b' \to S(C)$ is one of the four beachballs upstairs mapping to
it, each edge of $b$ maps to a path in $C$ of length $d-1$, which is covered by two paths in $S(C)$
of length $d-1$. Each such path is determined (given its image in $C$) by a word of length $d$
in the alphabet $0,1$. Our beachballs upstairs have the property that the pair of lifts of each
edge of $b$ are described by the {\em same} pair of $01$ words. Moreover, the collection of 4
pairs of $01$ words describing the lifts of $b$ are the same for all $b \to C$. So to describe
the lens move we just need to give 4 sets of pairs of $01$ words of length $d$. We call these
{\em height pairs}. An explicit formula for $d>2$ is given by
$$\text{ height pairs are } \begin{cases}
(00^{d-2}0, \quad 01^{d-2}0) \\
(0(01)^{\frac{d-2}{2}}1, \quad 0(10)^{\frac{d-2}{2}}1) \\
(1(1100)^{\frac{d-2}{4}}x, \quad 1(0011)^{\frac{d-2}{4}}x) \\
(1(1001)^{\frac{d-2}{4}}(1-x), \quad 1(0110)^{\frac{d-2}{4}}(1-x))
\end{cases}$$
where $x$ is $0$ or $1$ depending on the parity of $d$, and
where an expression like $w^{p/q}$ for $w$ a word of length $q$ means the initial string
of length $p$ of the word $w^\infty:=www\cdots$. The meaning of this formula is best 
explained by a picture; see Figure~\ref{height_graph}.

The case $d=2$ is degenerate; the case $d=3$ gives rise to the height pairs
$$(000,010), (001,011), (111,101), (110,100)$$
which reproduces four of the beachballs in the lens gluing for $d=3$ described above (the other
four beachballs are in the preimage of the other beachball in the hypercube gluing).

Note that there is an ambiguity in our choice of labels of each pair of vertices of $S(C)$
over $v$ by $v0$ and $v1$; thus the group $\Z/2\Z^{\text{vert}(C)}$ acts by automorphisms
of $S(C)$ over $C$. Thus although this is not evident in the formulae,
under this symmetry group all four height pairs are in the same orbit.

\begin{figure}[htpb]
\labellist
\small\hair 2pt
\endlabellist
\centering
\includegraphics[scale=0.6]{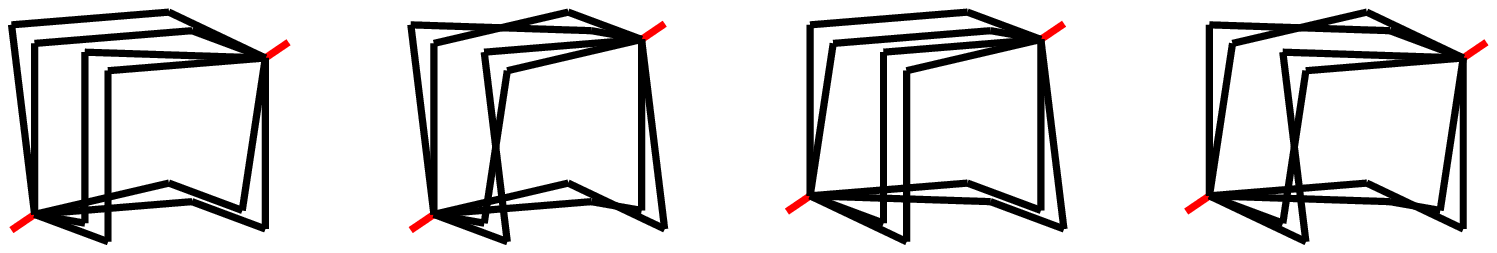}
\caption{Four beachballs in the $d=4$ lens gluing which project to one beachball in the
$d=3$ hypercube gluing.}\label{4lens}
\end{figure}

The case $d=4$ is hard to draw without the diagram becoming cluttered; Figure~\ref{4lens} shows
the 4 beachballs of degree 6 in the $d=4$ lens gluing which project to one beachball of degree 3 
in the $d=3$ hypercube gluing.

\subsection{Regularity}

The third condition for $(m,d)$-regularity is that each component of $L$ has exactly $m$ vertices
which map to high-valent vertices of $\Sigma$. The only move which adjusts the number of
such vertices on each component of $L$ is the elimination move. This move takes pieces
which may be drawn from any component $L$; so we simply take enough disjoint copies of $L$, and
symmetrize the components from which the moves are drawn, so that we apply this move the
same number of times to every component. This will ensure the third condition for $(m,d)$-regularity.

\subsection{Cocycle condition}

The cocycle is sensitive to more combinatorial data than we have been using so far; it depends not
only on the set of labels appearing as the edges in a beachball, but the way in which these
labels come from consecutive even segments on each of $\d$ blocks in a compatible $\d$-tuple.
Gluing different compatible $\d$-tuples gives rise to the same beachball collections but with
different cocycles (relative to a trivialization along the odd segments of the blocks). If the
original gluing was done randomly, there are pairs of compatible segments whose gluings give
rise to the same set of beachballs, but with cocycles differing by any element in the symmetric
group $S_{d-1}$ in the simplicial case, or any element of $BC_{d-1}$ in the cubical case; 
adjusting the gluing by interchanging elements of these pairs gives the same
collection of beachballs and the same remainder (and therefore the same gluing problem) but with
the cocycle adjusted in the desired away. Boundedly many moves of this kind for each component
produces an $(m,d)$-regular simplicial or cubical 
spine with all topological edges of length at least $\lambda$.

\subsection{Few generators}\label{subsection:simplifying}

We now indicate how to modify our arguments and constructions so that they hold without
our simplifying assumptions.

The first assumption --- that the length $n$ of the relators is divisible by some big fixed
constant $\lambda \cdot N$ (where also $\d$ divides $\lambda$) --- is easy to dispense with. We let $\rho$ be the remainder when
dividing $n$ by $\lambda \cdot N$, then find $\d$ copies of a subword of $r$ of length
$\lambda \cdot N + \rho$ so that it is legal to glue these $\d$ words in $\d$ distinct copies
of $r$; what is left after this step is a graph with edges whose lengths are divisible by
$\lambda \cdot N$, and the rest of the construction can go through as before.

The second assumption --- that the number $k$ of free generators satisfies $2k-1\ge 2d+1$ ---
can be finessed by looking more carefully at the construction. It can be verified that we do
not really use the hypothesis that the word $r$ is random with respect to the uniform
measure on reduced words of length $n$; instead we can make do with a much a weaker hypothesis,
namely that we generate random words by {\em any} stationary ergodic Markov
process, whose subword distribution of any fixed length has enough symmetry, and such
that there are at least $\d+1$ letters that may follow any legal substring. 

The relevance of this is as follows: a random reduced word of length $n$
in a free group of rank $k$ can be thought of as a random reduced word of length $n/s$
in a free group of rank $k(2k-1)^{s-1}$ (whose symmetric generating set is the
set of reduced words in $F_k$ of length $s$) generated by a certain stationary ergodic Markov process. We
can perform matching of subwords and build spines in this new ``alphabet''; the result will
not immerse in $X$, since there might be folding of subtrees of diameter at most $2s$.
But the spine will still have the crucial properties that the fundamental group of
the 2-complex $\overline{M}$ is commensurable with $\Delta(m,d)$, and that the number of
immersed paths of length $\nu$ grows like $(2k-1)^{s\nu/\lambda}$ where $\lambda$ is as
big as we like. Taking $s$ big enough to ``simulate'' a group of rank at least
$d+1$, we can then take $\lambda$ much bigger so that $s/\lambda$ is arbitrarily small.
This small exponential growth rate of subpaths is the key to proving that 
$\pi_1(\overline{M})$ maps injectively with quasiconvex image.

\subsection{Conclusion of the theorem}\label{subsection:conclusion}

The remainder of the argument is almost indentical to the arguments in \cite{Calegari_Walker} and
\cite{Calegari_Wilton}, and makes use in the same way of the mesoscopic small cancellation theory
for random groups developed by Ollivier \cite{Ollivier}. It depends on two ingredients:
\begin{enumerate}
\item{a {\em bead decomposition}, exactly analogous to that in \cite{Calegari_Walker}, \S~5 
or \cite{Calegari_Wilton}, \S~4; and}
\item{a {\em mesoscopic small cancellation argument}, which applies techniques of Ollivier
as in \cite{Calegari_Walker}, \S~6 or \cite{Calegari_Wilton}, \S~6.}
\end{enumerate}

We very briefly summarize the argument below just to indicate that no new ideas are needed, but
refer the reader to the cited references for details.

A bead decomposition is an $(m,d)$-regular simplicial or cubical spine $\Sigma$ which is made up of
$O(n^\delta)$ pieces of size $O(n^{1-\delta})$ arranged in a circle, with each piece separated from
the next by a {\em neck} (an unusually long glued segment) of length $C\log{n}$. We can construct a
spine as in \cite{Calegari_Wilton}, \S~4 as follows. Fix some small positive constant
$\delta$ and write $r$ as a product
$$r=r_1s_1r_2s_2\cdots r_ms_m$$
where each $r_i$ has length approximately $n^{1-\delta}$ and each $s_i$ has length approximately $n^\delta$, and
we choose lengths so that $m$ is divisible by $\d$. We
then fix a small positive constant $C<\delta/\log(2k-1)$ and look for a common subword $x$ of length
$C\log{n}$ in $s_i,s_{i+m/\d},\cdots,s_{i+(\d-1)m/\d}$ with indices taken mod $m$,
and then glue these distinct copies of $x$
into unusually long segments (called {\em lips}) of what will become the spine $\Sigma$.

The lips partition the remainder of $L$ into subsets $b_i$ each with mass $O(n^{1-\delta})$, and
we can extend the partial gluing along the lips into an $(m,d)$-regular simplicial or cubical spine
by gluing intervals in each $b_i$ independently. As in  \cite{Calegari_Walker}, \S~5 
and \cite{Calegari_Wilton}, \S~4 the resulting spine $\Sigma$ has the property that for any positive
$\beta$, any immersed segment $\gamma \to \Sigma$ with length $\beta n$ whose image in $X$
lifts to $r$ or $r^{-1}$, already lifts to $L$ (i.e.\/ it appears in the boundary of a disk of
$\overline{M}$). From this and the fact that a random 1-relator group is $C'(\mu)$ for any positive
$\mu$ it follows that $\overline{M} \to K_r$ is $\pi_1$-injective and its image is quasiconvex.

\medskip

Finally, the fact that the valence of $\Sigma$ is uniformly bounded and the length of every segment
is at least $\lambda$ where we may choose $\lambda$ as big as we like, implies that 
$\overline{M} \to K_r$ stays $\pi_1$-injective and quasiconvex when we attach the disks corresponding
to the remaining $(2k-1)^{nD}-1$ random relations. This depends on Ollivier's mesoscopic small 
cancellation theory \cite{Ollivier}; 
the argument is exactly the same as that in \cite{Calegari_Walker}, \S~6 and 
\cite{Calegari_Wilton}, \S~6.

This completes the proof of the Superideal Simplex Theorem and Superideal Cube Theorem.

\section{Acknowledgments}
Danny Calegari was supported by NSF grant DMS 1405466. The \LaTeX\/ macro for the symbol
$\Square$ was copied from {\tt pmtmacros20.sty} by F. J. Yndur\'ain. I would like to thank
the anonymous referee for their helpful comments.

\end{document}